\documentclass[11pt]{article}

\usepackage{a4wide,graphicx,amssymb,amsthm}
\usepackage{epsf}

\def\al{\alpha}
\def\be{\beta}

\def\de{\delta}

\def\om{\omega}

\def\map{\rightarrow}
\def\bq{\begin{equation}}
\def\eq{\end{equation}}

\def\rc{{\mathbb R}}
\def\zc{{\mathbb Z}}
\def\fg{{\mathfrak g}}

\def\ti{\tilde}
\def\Om{\Omega}

\newtheorem{thm}{Theorem}
\newtheorem{df}{Definition}

\def\zv{\bfseries\itshape}

\def\Hom{\mbox{\it Hom}}

\begin{document}

\begin{center}
{\Large \scshape Noncommutative differential forms and quantization of the odd symplectic category}
\vskip 7mm

Pavol \v Severa\\{\it Dept.~of Theoretical Physics\\Comenius
University\\Bratislava, Slovakia}\\{\tt
severa@sophia.dtp.fmph.uniba.sk} \vskip 4mm
\end{center}

\begin{abstract}
There is a simple and natural quantization of differential forms on odd Poisson supermanifolds,
given by the relation $[f,dg]=\{f,g\}$ for any two functions $f$ and $g$. We notice that this
non-commutative differential algebra has a geometrical realization as a convolution algebra of
the symplectic groupoid integrating the Poisson manifold.

This quantization is just a part of a
 quantization of the odd symplectic
category (where objects are odd symplectic supermanifolds and
morphisms are Lagrangian relations) in terms of $\zc_2$-graded
chain complexes. It is a straightforward consequence of the theory
of BV operator acting on semidensities, due to H.~Khudaverdian.

{\em MSC2000:} 58A10, 53D55, 53D05, 58A50

{\em Keywords:} differential forms, symplectic category, Batalin-Vilkovisky operator
\end{abstract}

\section{Introduction}

There is a well-known analogy (formulated by A.~Weinstein \cite{We}) between
symplectic manifolds and vector spaces:

\begin{center}
\begin{tabular}{r|l}
symplectic manifold & vector space\\
Lagrangian submanifold & vector\\
product & tensor product\\
opposite sympl.~form & dual space\\
Lagrangian relation & linear map
\end{tabular}
\end{center}

\noindent It is very fruitful, though it is just an analogy,
i.e.~there is no functor from the symplectic category to the
category of vector spaces.

However, there is such a quantization functor, if we substitute
symplectic manifolds with odd symplectic supermanifolds, and
vector spaces with $\zc_2$-graded chain complexes,
i.e.~$\zc_2$-graded vector spaces $V=V_0\oplus V_1$ with a
differential $D:V_0\map V_1$, $V_1\map V_0$, $D^2=0$. Namely, to
an odd symplectic supermanifold $Y$ it associates the vector space
of semidensities on $Y$, with $D$ the BV operator (due to
H.~Khudaverdian).
Before giving the details, let us look at the result of
quantization of odd symplectic groupoids.

\section{Differential forms on odd Poisson supermanifolds}

Let $X$ be a supermanifold and $\Om(X)$ the differential graded algebra of differential forms on $X$. By a {\em filtered deformation of $\Om(X)$} we mean a differential filtered algebra $A$ (i.e.~a differential algebra $A$ with an increasing filtration $F^0\subset F^1\subset\dots \subset A$ such that $dF^i\subset F^{i+1}$) with an isomorphism between $\Om(X)$ and the differential graded algebra $\mathit{Gr}A$ associated to $A$ (where $(\mathit{Gr}A)^i=F^i/F^{i-1})$. We have the following simple theorem:
\begin{thm}
There is a natural bijection between (isomorphism classes of) filtered deformations of $\Om(X)$ and odd Poisson structures on $X$. It is given by the commutation relation
\bq[f,dg]=\{f,g\},\eq where $f$ and $g$ are any functions on $X$, $[,]$ is the supercommutator in the deformed algebra and $\{,\}$ is the odd Poisson structure.
\end{thm}
\begin{proof}Notice that $F^0=C^\infty(X)$. Since $\Om(X)$ is generated
(as an algebra) by functions on $X$ and by their differentials, the same is true for $A$ (as is easily proved by induction). Finally, for any $f,g\in F^0$ we have $[f,dg]\in F^0$ (because $\Om(X)$ is graded commutative, i.e.~there $[f,dg]=0$). The algebra $A$ is thus known (up to canonical isomorphism) if we know the function $\{f,g\}:=[f,dg]$ for every two functions $f$ and $g$. It is straightforward to verify that $\{f,g\}$ has to be an odd Poisson structure. The converse (i.e.~given an odd Poisson structure, the formula (1) gives a fitered deformation of $\Om(X)$) can also be easily verified directly, but we give a more conceptual proof of this fact, using quantization of odd symplectic groupoids, later.
\end{proof}

If $\pi$ is an odd Poison structure on $X$, we denote the corresponding filtered deformation of $\Om(X)$ by $\Om_\pi(X)$.
Since our construction is natural, for any Poisson map $X_1\map
X_2$ we have a pullback map $\Om_{\pi_2}(X_2)\map\Om_{\pi_1}(X_1)$
preserving all the structure; hence, e.g., if $X$ is a Poisson Lie
group, $\Om_\pi(X)$ is a (differential filtered) Hopf algebra.

Let us give some examples of odd Poisson manifolds and of their
algebras $\Om_\pi(X)$. If $X=\Pi T^*M$ with its canonical odd
symplectic form,  $\Om_\pi(X)$ is the algebra of differential
operators acting on differential forms on $M$. The operator
corresponding to a function on $X$, i.e.~to a multivector field
$m$ on $M$, is $i_m$. The differential is the supercommutator with
the de Rham $d$.

As another example, if $\fg$ is a Lie algebra, take $X=\Pi\fg^*$
with its Kirillov-Kostant Poisson structure. Geometrically,
$\Om_\pi(X)$ is the convolution algebra of deRham currents on the
group $G$, supported at $1\in G$. More algebraically, it is the
crossed product of $\mathfrak A\fg$ with $\bigwedge\fg$; the
differential maps identically $\fg\subset\bigwedge\fg$ to
$\fg\subset\mathfrak A\fg$.

These two examples were special cases of $X=\Pi A^*$, where $A\map
M$ is a Lie algebroid. Lie theory is a source of other interesting
odd Poisson structures on graded supermanifolds. Recall these
important observations of A.~Vaintrob \cite{Va}: A Lie algebroid
structure on $A$ is equivalent to a degree $-1$ odd Poisson
structure $\pi$ on $A^*[1]$ (this is the one we just
 mentioned) and also to a degree $1$ odd vector field $Q$ on $A[1]$ with $Q^2=0$. A bialgebroid
 structure on $A$ is both $\pi$ and $Q$ on $A^*[1]$, such that ${\cal L}_Q\pi=0$. They can be
 combined to a single Poisson structure $\ti\pi=\pi+Q\partial_t$ on $A^*[1]\times\rc[2]$ (here
 $t$ is the coordinate on $\rc[2]$; recall that an odd Poisson structure on $X$ is an odd
 quadratic function $\pi$ on $T^*X$ such that $\{\pi,\pi\}=0$; $Q\partial_t$ is understood
 as a product of two linear functions on $T^*(\Pi A^*\times\rc)$).

A {\em quasi-bialgebroid} structure on $A$ is a graded principal $\rc[2]$-bundle $X\map A^*[1]$
with a $\rc[2]$-invariant odd Poisson structure $\ti\pi$ of degree $-1$. Choosing
a trivialization $X=A^*[1]\times\rc[2]$ we get a decomposition
$\ti\pi=\pi+Q\partial_t+\phi\partial_t^2$, where $\pi$ is a degree $-1$ odd Poisson structure
on $A^*[1]$, $Q$ a degree 1 odd vector field on $A^*[1]$ and $\phi$ a degree 3 odd function on
$A^*[1]$ ($\pi$, $Q$ and $\phi$ satisfy the obvious equations coming from $\{\ti\pi,\ti\pi\}=0$).
It is not clear to me what role the differential algebra $\Om_{\ti\pi}(X)$ may play e.g.~in the
problem of quantization of quasi-bialgebroids.

\section{Khudaverdian's BV operator}

In this section we recall several theorems of H.~Khudaverdian
\cite{Khu}.

Let $Y$ be an odd symplectic manifold and let $x^i$, $\xi_i$ be
local Darboux coordinates (i.e.~$\om=dx^i\wedge d\xi_i$). Let
$$\Delta={\partial^2\over\partial x^i\partial\xi_i},$$
understood as a differential operator from semidensities to
semidensities. Then
\begin{thm}
\begin{enumerate}
\item{$\Delta$ is odd and $\Delta^2=0$}

\item{$\Delta$ is formally selfadjoint}

\item{${\cal L}_{X_f}=[\Delta,f]$, where $X_f$ is the Hamiltonian
vector field generated by $f$ and $[,]$ the supercommutator}

\item{$\Delta$ is independent of the choice of Darboux
coordinates}
\end{enumerate}
\end{thm}
\begin{proof}{This is just a sketch: 1.~and 2.~are evident and 3.~can
be directly computed. From 1.~and 3.~we get that $[{\cal
L}_{X_f},\Delta]=0$, i.e. $\Delta$ is invariant under Hamiltonian
diffeomorphisms, hence we get 4.}
\end{proof}
Notice that Hamiltonian diffeomorphisms act trivially of the
cohomology of $\Delta$, since by 3.~they are homotopies. If $Y=\Pi
T^*M$ then one can identify differential forms on $M$ with
semidensities on $Y$ using Fourier transform along the fibres of
$\Pi TM$; then $\Delta$ becomes $d$ and 3.~contains as a special
case Cartan formula.

If $\al$ and $\be$ are semidensities on $Y$, we set
$$(\al,\be)=\int_Y\al\be$$
(provided the integral is well defined). Since $\Delta$ is
formally selfadjoint, if $\al$ is $\Delta$-closed and $\be$ is
$\Delta$-exact then $(\al, \be)=0$, i.e.~$(,)$ is well defined on
cohomology classes (provided the appropriate care is taken for the
finiteness of the integral).

Let $L\subset Y$ be a Lagrangian submanifold. We define a
$\de$-like generalized semidensity $\de_L$ supported on $L$: in
local Darboux coordinates $x^i$, $\xi_i$, with $L$ given by the
equations $\xi_i=0$, we set $\de_L=\prod_i\de(\xi_i)$.
\begin{thm}
$\de_L$ is $\Delta$-closed and independent of the choice of
coordinates
\end{thm}
\begin{proof} Closedness is obvious and independence can be
computed.\end{proof}

Since Hamiltonian flows generate homotopies on semidensities, we
know that $\de_L$ and $\de_{L'}$ lie in the same cohomology class
whenever $L$ and $L'$ can be connected by a Hamiltonian
diffeomorphism. Hence, if $\al$ is closed (and appropriate
integrals are finite) then $(\de_L,\al)=(\de_{L'},\al)$; this is the geometrical
basis for BV quantization.

\section{Quantization functor}

\begin{df}
The objects of {\zv the quantum odd symplectic category (QOSC)}
are odd symplectic manifolds and $\Hom(Y_1,Y_2)$ is the
$\zc_2$-graded complex of generalized semidensities on $\bar
Y_1\times Y_2$ (\/$\bar Y_1$ denotes $Y_1$ with the opposite
sympectic form); composition of morphisms is given by integration.
\end{df}

Notice that the composition is not always defined, hence we don't
really have a category. This is the same problem as with the
symplectic category; we'll keep it in mind and ignore it.

Composition is $\Delta$-equivariant since $\Delta$ is formally
selfadjoint. The quantization of an odd symplectic manifold $Y$
can be defined as $\Hom(pt,Y)$, i.e.~as the space of generalized
semidensities on $Y$.

\begin{df}
The {\zv quantization functor} from the odd symplectic category to
QOSC is the identity on objects and it maps any Lagrangian
relation $L$ to $\de_L$.
\end{df}

Of course, we have to prove that composition of $L$'s corresponds
to composition of $\de_L$'s. This is quite simple. If two relations $L_1$ and $L_2$ can be composed (i.e.~they satisfy the transversality condition) and their composition is $L$, the composition of $\de_{L_1}$ with $\de_{L_2}$ is a $\de$-like distribution supported at $L$, i.e.~a multiple of $\de_L$. Since both $\de_{L_1}$ and $\de_{L_2}$ are closed, so must be their composition, i.e.~it is a {\em constant} multiple of $\de_L$. Finally, we have to verify that the constant is 1. It is enough to do it at a single point. For example, we can locally deform $L_1$ and $L_2$ so that they are both given by setting some of chosen Darboux coordinates to 0. For such $L$'s the constant is clearly 1, therefore it is 1 also in the undeformed parts of $L$'s, and so it is 1 for the original $L$'s as well.

\section{Differential forms on odd
Poisson supermanifolds II}

Now the following dream is fulfilled: whenever we have an
algebraic structure in the odd symplectic world, we get a
corresponding structure in the world of $\zc_2$-graded chain
complexes. As the most obvious example, let us start with an odd
Poisson manifold $X$ and let $Y$ be its (local) odd symplectic
groupoid; the quantization of $Y$ will be a differential algebra.

The quantization of $Y$ is the space of generalized semidensities
on $Y$. We choose a reasonable subspace so that the product in the
algebra (given by composition with $\de_L\in\Hom(Y\times Y,Y)$,
where $L$ is the graph of the product in the groupoid $Y$) is well
defined. A cheap choice is to take smooth multiples of $\de_X$ and
of its derivatives (here $X\subset Y$ is the Lagrangian
submanifold of units of the groupoid $Y$). In this way we get a
differential algebra with an increasing filtration (given by the
degree of the distributions). As one easily sees, it is the
algebra $\Om_\pi(X)$.

\end{document}